# Householder Meets Student

John H. Elton and Andrew B. Gardner

June 2022

**Summary**.  Householder's algorithm for the QR factorization of a tall thin $n \times p$ full-rank matrix $X$ has the added bonus of producing, with no extra work, a matrix $M$ with orthonormal columns that are a basis for the orthocomplement of the column space of $X$ . We give a simple formula for $M^T x$ when $x$ is in that orthocomplement. The formula does not require computing $M$ , it only requires the R factor of a QR factorization. This is used to get a remarkably simple computable concrete representation of "independent residuals" in classical linear regression. For "Student's" problem, when $p = 1$, if $R_j = Y_j - \overline{Y}$, $j = 1, ..., n$ are the usual (non-independent) residuals, then $W_j = R_{j+1} - \dfrac{1}{\sqrt{n+1}} R_1, j = 1, ..., n-1$, gives $n-1$ i.i.d mean-zero normal variables whose sum of squares is the same as that of the $n$ residuals.  Those properties of this formula can (in hindsight) easily be verified directly, yielding a new simple and concrete proof of Student's theorem.  It also gives a simple way of generating $n-1$ exactly mean-zero i.i.d. samples from $n$ samples with unknown mean. Yiping Cheng[1] exhibited concrete linear combinations of the $Y_i's$ with these properties, in the context of a constructive proof of Student's theorem, but that representation is not so simple.

Analogous results are obtained for regression when there are more predictors, giving a very simple computable concrete formula for $n-p$ i.i.d "independent residuals" with the same sum of squares as that of the usual $n$ non-independent residuals.  A connection with Cochran's theorem is discussed.

jelton@bellsouth.net; elton@math.gatech.edu

## 1.  Introduction.

**Householder's  QR Factorization.**  Ask a student how to orthogonalize a set of vectors, and the answer will be, "Gram-Schmidt".  A good answer.  Ask how to get an orthonormal basis for the orthocomplement of the span of a set of vectors, say the columns of matrix $X$ , and the answer might be, "Find a basis for the null space of $X^T$ , and Gram-Schmidt that."  That, or something like that, is what you'll see on math.stackexchange.com.  It is correct in the sense that it will probably get you credit on a linear algebra exam, but it isn't a good answer, if you have a few vectors in a high-dimensional space, i.e., $X$ is "tall and thin", say $n \times p, p \ll n$ .  It's going to take you $O\left(n(n-p)^2\right) \sim O\left(n^3\right)$ operations.  But here is a conceptual geometric argument showing that's no good.  Consider isometrically transforming (e.g., rotating or reflecting) $\mathbb{R}^n$ so that $e_1$ , the first unit basis vector, lines up with $x_1$ , the first column of $X$ .  The rest of the frame comes along for the ride.  Next transform so that the new $e_2$ lines up with the part of $x_2$ that is perpendicular to $x_1$ , keeping the already-transformed $e_1$ in place, lined up with $x_1$. Continue this way until the first $p$ unit basis vectors of the natural frame have been transformed to span the columns of $X$. We're done! The other $n - p$ vectors of the natural frame have been transformed to be an orthonormal basis for $col(X)^\perp$. In matrix terms, if we actually want these vectors, we could multiply the $p$ transformation matrices together, and the last $n - p$ columns of the transpose of the result is our answer. If the transformations are reflections or rotations, this can be done in $O\left(n^2 p\right)$ time, because the matrices are of the form $I$ minus a rank one matrix for reflections, or minus a rank two matrix for rotations. Reflections might not seem as natural to our geometric intuition for moving the frame around, but they are more efficient.

That's cute, but this already appears implicitly (without our little physical model) in your linear algebra book that discusses the famous Householder $QR$ factorization ( e.g. Strang[**4**, p. 361-363]). It's just that this orthocomplement aspect of the algorithm doesn't show up in the exercises, so it's not as well-known. The Householder algorithm uses a product of $p$ reflections to obtain $U^T X = \begin{bmatrix} T \\ 0 \end{bmatrix}$ for an orthogonal matrix $U$ and an upper triangular $p \times p$ matrix $T$, non-singular because we assume $X$ is of full rank. Section 2 will review the steps of the algorithm to establish notation. There is a choice of sign at each step; one choice avoids catastrophic cancellation, and that is the choice usually made in practice (as done in e.g. [**5**, p. 73]). Refer to this as the *standard sign choice*; see section 2 for the definition. The relevance of that choice for this article is that it will guarantee a certain matrix of interest (in Theorem 1 below) is non-singular. Write $X = U \begin{bmatrix} T \\ 0 \end{bmatrix} = \begin{bmatrix} U_1 & U_2 \end{bmatrix} \begin{bmatrix} T \\ 0 \end{bmatrix} = U_1 T$, which exhibits $X$ as the product of a matrix with orthonormal columns and an upper triangular matrix, called a $QR$ factorization, and implies $U_1 = X T^{-1}$. Since $col(X) = col(U_1)$, the columns of $U_2$ are already an orthonormal basis for $col(X)^\perp$. That seems like a distinct advantage of using the Householder algorithm for orthonormalization: you get that orthonormal basis for the orthocomplement for free.

However, our interest here is not in having our hands on those columns of $U_2$, it is in computing $U_2^T x$ when $x \in col(X)^\perp$. $U_2$ is $n \times (n - p)$, so direct multiplication would take $O(n(n - p)) \sim O\left(n^2\right)$ time. Obviously that is wrong: One can compute $U^T x$ in $O(np)$ time, because $U$ is the product of $p$ reflections, and then just ignore the first $p$ rows to get $U_2^T x$. But Theorem 1 gives a simple formula for $U_2^T x$ that does not even use the reflections explicitly, it only uses $T$.

In fact, our interest is not really computation. Our interest in having a low-complexity formula will be for insight into the residuals in classical linear regression.

For a matrix (or column vector) $A$ with $n$ rows, let $A^{(p)}$ and $A_{(p)}$ denote, respectively, the first $p$ and last $n - p$ rows of $A$, so that $A = \begin{bmatrix} A^{(p)} \\ A_{(p)} \end{bmatrix}$, in block matrix notation.

**Theorem 1.** *Let* $X$ *be* $n \times p$ *with* $p < n$ , *and of full rank.* *Let* $U^T X = \begin{bmatrix} U_1 & U_2 \end{bmatrix}^T X = \begin{bmatrix} T \\ 0 \end{bmatrix}$ *be the result of the Householder algorithm for QR factorization, where* $U$ *is orthogonal,* $U_1$ *is* $n \times p$ , $T$ *is* $p \times p$ *and upper triangular, and the standard sign choice is made at each step of the algorithm.* *Then* $T - X^{(p)}$ *is non-singular, and for* $x \in col(X)^{\perp}$ ,

$$U_2^T x = x_{(p)} + X_{(p)} \left( T - X^{(p)} \right)^{-1} x^{(p)}$$

If arbitrary sign choices are made in the Householder algorithm, it is possible that $T - X^{(p)}$ is singular, but in section 2 it is shown that there is nevertheless a $p \times p$ matrix $S$ that is computable in $O\left( p^3 \right)$ time, such that $U_2^T x = x_{(p)} + X_{(p)} S x^{(p)}$ . This $S$ has the same rank as $T - X^{(p)}$ , and is a type of generalized inverse of $T - X^{(p)}$ .

Notice that the formula does not require actually having your hands on the matrix $U_2$ . One can compute the action of $U_2^T$ on an element of $col(X)^{\perp}$ without having $U_2$ ; it is lurking in the background, but not actually needed. Theorem 2 does not even mention Householder reflections, though they will be implicitly there, as the proofs in section 2 will show.

**Theorem 2.** *Let* $X$ *be* $n \times p$ *with* $p < n$ , *and of full rank.* *Let* $C$ *be a non-singular* $p \times p$ *matrix such that* $XC^{-1}$ *has orthonormal columns.* *Then there exists an* $n \times (n-p)$ *matrix* $M$ *whose columns are an orthonormal basis for* $col(X)^{\perp}$ , *and a* $p \times p$ *matrix* $S$ *having the same rank as* $C - X^{(p)}$ , *such that for all* $x \in col(X)^{\perp}$ ,

$$M^T x = x_{(p)} + X_{(p)} S x^{(p)}.$$

*If* $C - X^{(p)}$ *is non-singular, then* $S = \left( C - X^{(p)} \right)^{-1}$ *; in any case,* $S$ *can be computed in* $O(p^3)$ *time.*

*Remark.* Actually, it is always possible to avoid the singular case by adjusting $C$ . If $C - X^{(p)}$ is singular, there exists a diagonal matrix $D$ with $\pm 1$'s on the diagonal (i.e., $D^2 = I$ ) such that $DC - X^{(p)}$ is non-singular (see section 2). But $X(DC)^{-1} = XC^{-1}D^{-1}$ also has orthonormal columns, so

$$M^T x = x_{(p)} + X_{(p)} \left( DC - X^{(p)} \right)^{-1} x^{(p)}$$

for some $M$ having the properties described above. $D$ can be found in $O(p^3)$ time.

In the proof, $M$ will be the last $n - p$ columns of a product of $p$ Householder reflections, used in the Householder algorithm applied to $XC^{-1}$ . But the theorem may be applied without knowing anything about that. $C$ might have come from Gram-Schmidt, say.

There is nothing special about separating into the first $p$ and last $n - p$ rows: that is for notational convenience. Any $p$ rows may be used (of course, $T$ or $C$ would change accordingly).

**Statistics with Householder.** As far as we know, A.S. Householder never met W. S. Gosset (a.k.a. "Student", to anyone who has taken a statistics course). But we will let Householder's algorithm for QR factorization meet up with Student's problem, to give an interesting representation. The classical "Student's" theorem shows that for $Y_1, \ldots Y_n$ independent random variables with common distribution $N(\mu, 1)$, $\sum_{j=1}^n \left(Y_j - \overline{Y}\right)^2$ has the chi-square distribution with $n-1$ degrees of freedom. This appears in every first course in statistics for students having a math background, and it is the basis for the t-test that is taught to all students taking any kind of statistics course. A proof, however, is not obvious. Some books skip it, and probably many students at that level find the proofs that are given a little hard to grasp (if they attempt it); they are all abstract, in the sense that they do not actually exhibit a concrete, computable vector of $n-1$ i.i.d mean-zero normal variables whose sum of squares matches that of the residuals, $\sum_{j=1}^n \left(Y_j - \overline{Y}\right)^2$ .

Yiping Cheng [**1**] had the nice idea that for pedagogical reasons, it would be good to actually exhibit such a representation. He gives a formula for an $n \times (n-1)$ matrix $M$ with orthonormal columns which are orthogonal to $\mathbf{1}$ (the vector of ones), such that if $W = M^T \left(Y - \mathbf{1}\mu\right)$, then $\sum_{j=1}^{n-1} W_j^2 = \sum_{j=1}^n \left(Y_j - \overline{Y}\right)^2$ . Because of the orthonormal columns, $W \sim N(0, \sigma^2 I_{n-1})$, where $I_k$ denotes the $k \times k$ identity matrix. Since $M^T \mathbf{1} = 0$, actually $W = M^T Y$ which does not involve the unknown $\mu$ , so we have a computable representation: $W$ is a *statistic*. In section 3 we will show a simple way to arrive his matrix in the context of our discussion of Cochran's theorem. But the calculation of $M^T Y$ using that matrix takes $O\left(n^2\right)$ operations, so is not particularly simple.

The theorems above, when $p = 1$, imply a surprisingly simple formula for a concrete representation in Student's problem. And for classical linear regression with more parameters, a simple concrete representation is also obtained. Here's the regression setup (see e.g. [2, p. 522]).

***Linear Regression Setup:*** $X$ is an $n \times p$ matrix with $p < n$, $\beta \in \mathbb{R}^p$ are unknown parameters, and $Y$ is a multivariate normal random vector with mean $X\beta$ and covariance matrix $\sigma^2 I_n$; that is, $Y - X\beta \sim N(0, \sigma^2 I_n)$ . Assume $rank(X) = p$, so $X$ is of full rank, with linearly independent columns. The normal equations for the maximum likelihood estimator $\hat{\beta}$ of the parameters are $X^T Y = X^T X \hat{\beta}$ , so $X^T \left(Y - X\hat{\beta}\right) = 0$ . Thus the columns of $X$ are perpendicular to the residuals: $X^T R = 0$ , where $R = Y - X\hat{\beta}$ are the residuals; that is, $R \in col(X)^\perp$ . The estimator $\hat{\beta} = \left(X^T X\right)^{-1} X^T Y$ is a linear function of $Y$ .

Theorem 3 follows from Theorems 1 and 2 above (see section 3).

**Theorem 3.** *Assume the linear regression model above. Let* $v = \left(T - X^{(p)}\right)^{-1} R^{(p)}$, *where* $T$ *is the upper triangular matrix in the Householder QR factorization of* $X$ , *using the standard sign choices. Or, if* $XC^{-1}$ *has orthonormal columns, let* $v = \left(C - X^{(p)}\right)^{-1} R^{(p)}$ *if* $C - X^{(p)}$ *is non-singular, or more* generally, $v = SR^{(p)}$ where $S$ is as in Theorem 2.  Then $W = R_{(p)} + X_{(p)}v \sim N(0, \sigma^2 I_{n-p})$, and $W^T W = R^T R$ .

As an alternate interpretation, let $\beta^* = \hat{\beta} - v$ , so $W = Y_{(p)} - X_{(p)}\beta^*$ .  The alternate interpretation views it as perturbing the usual least squares estimator to get $n - p$ "independent residuals".  Again, there is nothing special about using the first $p$ rows.

Theorem 3, with $p = 1$, may be applied to Student's problem, where $Y - X\mu \sim N\left(0, \sigma^2 I_n\right)$ with $X = 1$, and $\hat{\mu} = \overline{Y}$, and $R = Y - 1\overline{Y}$ .  Since $T$ and $C$ are $1 \times 1$ matrices, there is no difference between the versions, and clearly $C = -\sqrt{n}$ or $C = \sqrt{n}$ makes $XC^{-1}$ have an orthonormal column.  Since $X^{(1)} = 1$, $\left(C - X^{(1)}\right)^{-1} = -\dfrac{1}{\sqrt{n}+1}$ and $\left(C - X^{(1)}\right)^{-1} R^{(1)} = -\dfrac{1}{\sqrt{n}+1} R_1$ for the first case, so $W = R_{(1)} - 1\dfrac{1}{\sqrt{n}+1}R_1 \sim N(0, \sigma^2 I_{n-1})$ . The other sign gives $W = R_{(1)} + 1\dfrac{1}{\sqrt{n}-1}R_1$ .

It looks simpler just writing components:

**Corollary 4**.  *In Student's problem, let* $R_j = Y_j - \overline{Y}$, *and let*

$$W_j = R_{j+1} - \frac{1}{\sqrt{n}+1}R_1, \; j = 1,...,n-1 \; ; \text{ or let } W_j = R_{j+1} + \frac{1}{\sqrt{n}-1}R_1, \; j = 1,...,n-1 \; .$$

*Then the* $W_j$ *are i.i.d with* $W_j \sim N(0, \sigma^2)$, *and* $\sum_{j=1}^{n-1}W_j^{\,2} = \sum_{j=1}^{n}R_j^{\,2}$ .

Rather than interpreting this as a perturbation of $n-1$ of the residuals by the other one, it can be interpreted as a modification of the estimator of the mean.  Let $\mu^* = \hat{\mu} + \dfrac{1}{\sqrt{n}+1}R_1$, or alternatively,

$\mu^* = \hat{\mu} - \dfrac{1}{\sqrt{n}-1}R_1$ .  Then $W_j = Y_{j+1} - \mu^*, \; j = 1,...,n-1$ .

Since the $W_j$ are so simple to compute – essentially as easy as the usual residuals -  perhaps they are of interest interpreted as $n-1$ "independent residuals", rather than the usual $n$ non-independent residuals.  They have the same variance as the source.

This also could be looked at as a simple way of generating $n-1$ *exactly* mean zero i.i.d. normals from an i.i.d.  sample of size $n$ from a normal distribution with unknown mean.  Perhaps you have a source known to be exactly normal, but with unknown mean. Using this formula, you get mean zero samples (with variance the same as the source) without having to know the source mean, at the cost of only giving up one sample. Now, one could get $n/2$ independent mean zero samples by looking at $(Y_1 - Y_2)/\sqrt{2}$, $(Y_3 - Y_4)/\sqrt{2}, ..., (Y_{n-1} - Y_n)/\sqrt{2}$ with even $n$, but that wastes half the samples.

There is yet another way to write the result. Let $R_j^* = Y_j - \overline{Y_{(1)}}$, where $\overline{Y_{(1)}}$ is the average of the last $n-1$ variables; that is, $\overline{Y_{(1)}} = \frac{1}{n-1}\sum_{i=2}^n Y_i$. A little algebra shows that

$$R_{j+1} + \frac{1}{\sqrt{n-1}}R_1 = R_{j+1}^* + \frac{1}{\sqrt{n}}R_1^*, \quad j = 1, ..., n-1.$$ The interesting thing here is that $R_1^*$ is statistically independent from $R_{j+1}^*$, as may be easily shown. Since there is nothing special about which sample is used to "correct" the others, consider this model, where the order is reversed. Draw $Y_1, ..., Y_{n-1}$, and form the sample average $\overline{Y^{(n-1)}} = \frac{1}{n-1}\sum_{i=1}^{n-1} Y_i$ of those, and subtract that from $Y_1, ..., Y_{n-1}$, forming the ordinary residuals $R_j^* = Y_j - \overline{Y^{(n-1)}}$, $j = 1, ..., n-1$ for a sample of size $n-1$. Then draw one more sample, and correct $R_1^*, ..., R_{n-1}^*$ by adding an independent term, $\frac{1}{\sqrt{n}}R_n^*$, to each, which shatters their correlation.

Also for regression with $p > 1$, one could consider the $n-p$ $W_j$'s as truly independent residuals, having the same variance as the source. They are very simple to compute. From Theorem 3, written in coordinates, $W_j = R_{j+p} + \sum_{i=1}^p v_i R_i$, $j = 1, ..., n-p$. Or, this can be considered a way of generating an exactly mean zero i.i.d. normal sample of size $n-p$ from a sample of size $n$ taken from a source known to follow a linear regression model, but with unknown parameters, which are not needed in the computation.

Just as for Student's problem, to give another concrete example, the result may be written out in coordinates for the case of ordinary univariate (that is, intercept and slope) regression. The result is a very simple linear perturbation of the ordinary residuals. The notation is easier if it is set this up with orthonormal columns. Assume WLOG that the predictor $t$ satisfies $\sum_{j=1}^n t_j = 0$, $\sum_{j=1}^n t_j^2 = 1$, and with i.i.d. random variables $Y_j$ such that $Y_j - a - bt_j \sim N(0, \sigma^2)$, $j = 1, ..., n$, where $a, b$ are unknown parameters, and $R_j = Y_j - \hat{a} - \hat{b}t_j$, $j = 1, ..., n$. This is the regression case $p = 2$, with $X = \begin{bmatrix} \frac{1}{\sqrt{n}}\mathbf{1} & t \end{bmatrix}$ and $\beta = \begin{bmatrix} a\sqrt{n} & b \end{bmatrix}^T$, with $X$ already having orthonormal columns. Use Theorem 3 with $C = I$, which has no messing with signs but has a possible singular case; or, use $T$ coming from the Householder algorithm with the standard sign choices, which has no singular case. See section 3 for the details.

**Corollary 5**. *Assume the univariate regression setup above. Let*

$$W_j = R_{j+2} + (AR_1 + BR_2) + (CR_1 + DR_2)t_{j+2}, j = 1, ..., n-2,$$

*where either:*

*(a)* $\begin{bmatrix} A & B \\ C & D \end{bmatrix} = \dfrac{1}{(\sqrt{n}-1)(1-t_2)-t_1} \begin{bmatrix} 1-t_2 & t_1 \\ 1 & \sqrt{n}-1 \end{bmatrix}$ *if* $(\sqrt{n}-1)(1-t_2)-t_1 \neq 0$,

$\quad\begin{bmatrix} A & B \\ C & D \end{bmatrix} = \dfrac{1}{\sqrt{n}-1} \begin{bmatrix} 1 & 0 \\ 0 & 0 \end{bmatrix}$ *if* $(\sqrt{n}-1)(1-t_2)-t_1 = 0$; *or*

*(b)* $\begin{bmatrix} A & B \\ C & D \end{bmatrix} = \dfrac{-s}{(\sqrt{n}+1)+\left|(\sqrt{n}+1)t_2-t_1\right|} \begin{bmatrix} s+t_2 & -t_1 \\ -1 & \sqrt{n}+1 \end{bmatrix}$, *where* $s = \mathrm{sgn}\left((\sqrt{n}+1)t_2-t_1\right)$.

*Then the* $W_j$ *are i.i.d with* $W_j \sim N(0,\sigma^2)$, *and* $\sum_{j=1}^{n-2} W_j^{\,2} = \sum_{j=1}^{n} R_j^{\,2}$.

Defining $a^* = \hat{a} - AR_1 - BR_2, b^* = \hat{b} - CR_1 - DR_2$, this can be written as $W_j = Y_{j+2} - a^* - b^* t_{j+2}, j = 1,...,n-2$, interpreting the result as a modification of the usual parameter estimates.

The singular case in (a) can occur for any n, but in essentially only one way:

$t_1 = \dfrac{1}{\sqrt{n}}, t_2 = 1 - \dfrac{1}{\sqrt{n}(\sqrt{n}-1)}, t_j = -\dfrac{1}{\sqrt{n}(\sqrt{n}-1)}, j = 3,...,n$. It illustrates the case of a rank-one $S$ in

Theorem 2, with the solution described in Lemma 10 below. Numerically, the solution (b) coming from Householder with standard sign choices is better, because the solution in (a) blows up near the singular case.

**Statistics without Householder.** Now that we've tried to convince you that the Householder algorithm leads to wonderfully simple formulas for computable representations of the normal vectors that occur in regression theory, we observe that, in hindsight, one could have conceivably arrived at these formulas with just statistical reasoning.

For the Student problem, one does not even need to mention a matrix. It is conceivable that without any motivation from Householder, one would think to try for a solution of the form $W_j = R_{j+1} + cR_1, \ j = 1,...,n-1$, to satisfy the independence and sum of squares conditions. With just simple algebra, one may prove (see section 4):

**Theorem 6**. *In Student's problem, let* $W_j = R_{j+1} + cR_1, \ j = 1,...,n-1$, *where* $R_j = Y_j - \overline{Y}$. *Imposing the conditions that the* $W_j$ *be i.i.d with* $W_j \sim N(0,\sigma^2)$, *and* $\sum_{j=1}^{n-1} W_j^{\,2} = \sum_{j=1}^{n} R_j^{\,2}$, *there are exactly two solutions for* $c$: $\quad c = \dfrac{1}{\sqrt{n}-1} \ or \ \dfrac{-1}{\sqrt{n}+1}$.

This is surely the most elementary proof of Student's theorem, and it could be given as an exercise. But it lacks motivation. Perhaps one could have come up with this form by thinking this way: the $R_j$'s are slightly correlated, and if the same little amount of random variable $R_1$ is added to each of the one of the $R_j$'s for $j > 1$, they remain identically distributed, so just find the amount to add to bring their correlations to zero. That much is easy. By luck?, adding that bit also brings the sum of the $n-1$ squares up to match the sum of the $n$ squares of the $R_j$'s. Well, we didn't think of the form without having first done the Householder reflection method.

Looking at the matrix form of this solution does not seem to give any additional insight. $W = \begin{bmatrix} c\mathbf{1} & I_{n-1} \end{bmatrix} R$, with $\begin{bmatrix} c\mathbf{1} & I_{n-1} \end{bmatrix}$ in block matrix form. This is a very simple matrix, but it does not have orthonormal rows. Theory leads one to express $W = M^T R$ where $M$ has orthonormal columns which are perpendicular to $\mathbf{1}$, so that the properties of $W$ are apparent. Then it was found that $M^T R = \begin{bmatrix} c\mathbf{1} & I_{n-1} \end{bmatrix} R$ for two values of $c$, so $W$ can be computed in this simpler way, and $M$ is not needed. But starting with the form $W = \begin{bmatrix} c\mathbf{1} & I_{n-1} \end{bmatrix} R$ and finding $c$ as described in section 4, there is no apparent connection with the $M$ that came from Householder.

For the regression problem with $p \geq 2$, the same thing can be done: Assume a form for a potential answer, and then just find the solutions. This uses matrix notation, but nothing about Householder; just ask for the required independence and sum of squares.

**Theorem 7.** *Consider the problem of finding a $p \times p$ matrix $S$ such that for $W = R_{(p)} + X_{(p)} S R^{(p)}$, $W \sim N(0, \sigma^2 I_{n-p})$ and $W^T W = R^T R$. If $C$ is a $p \times p$ matrix such that $XC^{-1}$ has orthonormal columns and $C - X^{(p)}$ is non-singular, then $S = \left( C - X^{(p)} \right)^{-1}$ is a solution, and all non-singular solutions are of this form, for some such $C$.*

The proof, in section 4, does not use the Householder method or any other matrix theory, beyond just inverse matrices. Such a $C$ can be found by elementary means. This gives a proof, rather different from the usual one, of the classical result that $R^T R / \sigma^2$ has the $\chi^2(n-p)$ distribution.

I don't think we would have thought of this form without having first done the Householder algorithm analysis, though. And although it doesn't show up explicitly, a $W$ having this simple form must actually be $U_2^T R$ where $U_2$ is the last $n-p$ columns of a product of $p$ Householder reflections. The solution here in matrix form, showing it as a linear transformation of $R$, is $W = \begin{bmatrix} X_{(p)} S & I_{n-p} \end{bmatrix} R$, which is very simple but whose motivation does not seem transparent.

## 2. Householder's QR Factorization.

Let $X$ be an $n \times p$ matrix of full rank, with $p \leq n$. The case of interest here is when $p$ is small compared to $n$ (a tall thin matrix).

We will now describe the Householder algorithm for obtaining an orthonormal basis for the column space of $X$, establishing a notation where each successive step consists of applying an $n \times n$ Householder reflection matrix to the previous step. Householder's factorization algorithm is nowadays discussed in beginning linear algebra texts, thanks to Strang ( [**4**, p. 361]). The usual description has the reflection matrices and vectors decreasing in size at each step, with previously done columns unchanged at each step and not part of the notation. But for notational purposes it is easier to use $n \times n$ matrices and $n$-vectors at each step; the vector at step $k$ will be zero in components above the $k^{th}$. In numerical implementation the stored vectors could decrease in size at each step. It hardly matters for the tall thin case.

A Householder transformation, or elementary reflector, has matrix $H = I_n - \dfrac{2}{\|v\|^2} vv^T$ where $I_n$ is the $n \times n$ identity matrix and the vectors are $n$-vectors. It is a symmetric orthogonal matrix, called a *Householder matrix*. Given vectors $x$, $y$ of the same length, one can choose $v$ so that $Hx = y$: let $v = x - y$, "from" minus "to"; any non-zero scalar multiple gives the same $H$. It is easy to see that if $x \neq y$, $v$ has to be a scalar multiple of $x - y$, so that $H$ is uniquely determined by $x$ and $y$ if they are not the same. A picture makes everything clear.

For the Householder algorithm, the elementary reflectors are used to transform column vectors to have zeroes below the diagonal. Recursively define a sequence $H_1, ..., H_p$ of such matrices

$H_k = I_n - \dfrac{2}{\|v_k\|^2} v_k v_k^T$, with $v_k$ being zero in its first $k-1$ components, in such a way that

$H_k ... H_2 H_1 \begin{bmatrix} x_1 & x_2 & ... & x_k \end{bmatrix} = \begin{bmatrix} t_1 & t_2 & ... & t_k \end{bmatrix}$, $k = 1, ..., p$, where vector $t_i$ is zero below the $i^{th}$

component; $x_i$ is the $i^{th}$ column of $X$. Notice that

$H_{k+1} H_k ... H_2 H_1 \begin{bmatrix} x_1 & x_2 & ... & x_k \end{bmatrix} = H_{k+1} \begin{bmatrix} t_1 & t_2 & ... & t_k \end{bmatrix} = \begin{bmatrix} t_1 & t_2 & ... & t_k \end{bmatrix}$, since $v_{k+1}$ is zero in its

first k components, so $H_{k+1}$ behaves as the identity on these vectors. In other words, the later reflection operators do not change the triangular structure that has already been obtained in earlier steps. At the end, $H_p ... H_2 H_1 X = \begin{bmatrix} T \\ 0 \end{bmatrix}$, with $T = \begin{bmatrix} t_1 & ... & t_p \end{bmatrix}^{(p)}$, a $p \times p$ upper triangular matrix.

Since the Householder matrices are symmetric and orthogonal, $X = H_1 H_2 ... H_p \begin{bmatrix} T \\ 0 \end{bmatrix} = QT$, where $Q$ is the first $p$ columns of the orthogonal matrix $H_1 H_2 ... H_p$. That exhibits a so-called "QR" factorization of $X$ (using $T$ for the upper triangular matrix rather than $R$). The last $n - p$ columns of $H_1 H_2 ... H_p$ are an orthonormal basis for $col(X)^\perp$ (a nice feature that motivated the main results of this article).

Since $v_{k+1}$ is zero in its first k components, one could write $H_{k+1}$ in block matrix form as

$H_{k+1} = \begin{bmatrix} I_k & 0 \\ 0 & I_{n-k} - \dfrac{2 w_{k+1} w_{k+1}^T}{\|w_{k+1}\|^2} \end{bmatrix}$, where $w_{k+1} \in \mathbb{R}^{n-k}$, and $v_{k+1} = \begin{bmatrix} 0 \\ w_{k+1} \end{bmatrix}$. But the block matrix

notation would be clumsy for some of the later discussion, so using $v_{k+1}$ rather than $w_{k+1}$ is preferred.

Look at the steps of the algorithm. Start by choosing $H_1$ to map $x_1$, the first column of $X$, to a multiple of $e_1$, so $v_1 = x_1 \pm \|x_1\| e_1$ (the "to" minus "from" rule), and either sign is allowed to be chosen, but a certain sign choice is usually made in implementations for numerical reasons. If $v_1 = x_1 \pm \|x_1\| e_1$, then $H_1 x_1 = \mp \|x_1\| e_1$, so the first column is taken care of. Write $v_1 = x_1 + d_1 \|x_1\| e_1$ where $d_1 = \pm 1$.

Most implementations choose the sign that prevents a loss of significant bits by cancellation in subtraction, and take $d_1 = \text{sgn}(x_{11})$; call this the *standard sign choice*. This also guarantees that $v_1 \neq 0$

, since $x_1$ is not zero (the columns of $X$ are linearly independent). But $v_1$ could be zero for arbitrary choice of sign; in that case, $H_1$ would be the identity matrix. Let $G_0 = I_n$.

Now assume step $k$ is finished, with $G_k = H_k ... H_1$, such that $G_k \begin{bmatrix} x_1 & ,..., & x_k \end{bmatrix} = \begin{bmatrix} t_1 & ... & t_k \end{bmatrix}$ with $(t_i)_j = 0, j > i$ so that the result is upper triangular so far. The goal is to define $H_{k+1}$ so that

$$G_{k+1} x_{k+1} = H_{k+1} G_k x_{k+1} = \begin{bmatrix} (G_k x_{k+1})^{(k)} \\ c \\ 0 \end{bmatrix} = t_{k+1},$$ so that the first $k$ components of $G_k x_{k+1}$ are unchanged

and the result $t_{k+1}$ is zero below the $(k+1)^{st}$ component, and the triangular structure is continued.

$H_{k+1}$ is an isometry, which forces $|c| = \left\| (G_k x_{k+1})_{(k)} \right\|$. Using "from" minus "to", let

$$v_{k+1} = \begin{bmatrix} 0 \\ (G_k x_{k+1})_{(k)} \end{bmatrix} + d_{k+1} \left\| (G_k x_{k+1})_{(k)} \right\| e_{k+1} \text{ with } d_{k+1} = \pm 1, \text{ where the 0 in the first vector has } k$$

components; thus the first $k$ components of $v_{k+1}$ are zero, which implies $v_{k+1}^T t_i = 0, i \le k$ , and

$H_{k+1} t_i = t_i, i \le k$ : the previous columns are unchanged. So

$$G_{k+1} \begin{bmatrix} x_1 & ,..., & x_{k+1} \end{bmatrix} = H_{k+1} G_k \begin{bmatrix} x_1 & ,..., & x_{k+1} \end{bmatrix} = \begin{bmatrix} t_1 & ... & t_{k+1} \end{bmatrix} .$$

Note that $(G_k x_{k+1})_{(k)}$ cannot be zero, because if it were, $G_k x_{k+1}$ would be in the span of the columns of $G_k \begin{bmatrix} x_1 & ,..., & x_k \end{bmatrix} = \begin{bmatrix} t_1 & ... & t_k \end{bmatrix}$, but the columns of $X$ are assumed linearly independent so this is impossible. Note $v_{k+1}$ is not zero if one takes $d_{k+1} = \text{sgn} (G_k x_{k+1})_{k+1}$ (the standard choice); it would be zero if and only if the opposite sign is used and if also $\left\| (G_k x_{k+1})_{k+1} \right\| = \left\| (G_k x_{k+1})_{(k)} \right\|$, so that all the mass is located on component $k+1$. By induction, none of the $v_k's$ are zero if the standard sign choice is made. Also observe (used in what follows) that $(v_{k+1})_{k+1} = 0 \Leftrightarrow v_{k+1} = 0$ since cancellation of component $k+1$ in calculating $v_{k+1}$ implies all the mass of $(G_k x_{k+1})_{(k)}$ is on that component.

This completes the summary of the Householder factorization algorithm, using a notation where $H_k = I_n - \dfrac{2}{\|v_k\|^2} v_k v_k^T$ (or just $I_n$ if $v_k = 0$ ) has all matrices and vectors of size $n$ .

The following lemma will be needed in the proof of Theorem 1.

**Lemma 8.** *Let* $X = G_p^T \begin{bmatrix} T \\ 0 \end{bmatrix} = H_1 H_2 ... H_p \begin{bmatrix} T \\ 0 \end{bmatrix}$ *be the result of the Householder factorization algorithm*

*applied to* $X$ , *as just described, with* $H_k = I_n - \dfrac{2}{\|v_k\|^2} v_k v_k^T$ *(or* $I_n$ *if* $v_k = 0$ *) and* $T$ *upper triangular.*

*Then*

$$rank \left( T - X^{(p)} \right) = \# \left\{ 1 \le i \le p : v_i \ne 0 \right\}$$

*In particular, $T - X^{(p)}$ is non-singular if the standard sign choice is made at each step in the algorithm.*

*Proof.* To simplify the notation, assume that the non-zero $v_i$'s have been normalized.

$H_1 H_2 ... H_p = \prod_{i=1}^{p} \left( I_n - 2 v_i v_i^T \right) = I_n - \sum_{i=1}^{p} 2 v_i v_i^T \prod_{j=i+1}^{p} \left( I_n - 2 v_j v_j^T \right)$, so

$I_n - H_1 H_2 ... H_p = \sum_{i=1}^{p} v_i r_i^T$ where $r_i^T = 2 v_i^T \prod_{j=i+1}^{p} \left( I_n - 2 v_j v_j^T \right)$. Since $(v_j)_i = 0$ for $j > i$,

$(r_i)_i = 2(v_i)_i$, so $(v_i)_i \neq 0 \Rightarrow (r_i)_i \neq 0$. If the product in the definition of $r_i^T$ were expanded out, each term would end with a row vector $v_j^T$ for some $j \geq i$, and $(v_j)_k = 0$ for $k < j$, so $(r_i)_k = 0$ for $k < i$.

We have $I_p - \left( \text{first } p \text{ rows and cols of } H_1 H_2 ... H_p \right) = \sum_{i=1}^{p} v_i^{(p)} \left( r_i^{(p)} \right)^T$. But

$$ X = H_1 H_2 ... H_p \begin{bmatrix} T \\ 0 \end{bmatrix} \Rightarrow X T^{-1} = H_1 H_2 ... H_p \begin{bmatrix} I_p \\ 0 \end{bmatrix} \Rightarrow X^{(p)} T^{-1} = \text{first } p \text{ rows and cols of } H_1 H_2 ... H_p $$

So $I_p - X^{(p)} T^{-1} = \sum_{i=1}^{p} v_i^{(p)} \left( r_i^{(p)} \right)^T$. Note that $v_i^{(p)}$ and $r_i^{(p)}$ satisfy the conditions inherited from $v_i$ and $r_i$ : $(v_i^{(p)})_j = (r_i^{(p)})_j = 0, j < i$; $(v_i^{(p)})_i \neq 0 \Rightarrow (r_i^{(p)})_i \neq 0$; and $\left( v_i^{(p)} \right)_i = 0 \Rightarrow v_i^{(p)} = 0$. By the lemma which follows, $rank \left( I_p - X^{(p)} T^{-1} \right) = \# \{ 1 \leq i \leq p : v_i \neq 0 \}$. But the rank is not changed multiplying by the non-singular matrix $T$, so the theorem is proved. □

**Lemma 9.** *Let $c_i \in \mathbb{R}^p, r_i \in \mathbb{R}^p$ such that (i) $(c_i)_j = (r_i)_j = 0, j < i$; (ii) $(c_i)_i \neq 0 \Rightarrow (r_i)_i \neq 0, 1 \leq i \leq p$, and (iii) $(c_i)_i = 0 \Rightarrow c_i = 0, 1 \leq i \leq p$. **Then** $rank \left( \sum_{i=1}^{p} c_i r_i^T \right) = \# \{ i : c_i \neq 0 \}$.*

The proof, which is left to the reader, is by (backward) induction using elementary row operations.

Next comes the key result for this paper about the Householder algorithm, a formula representation of $H_p ... H_2 H_1$, which will lead to Theorems 1 and 2. To facilitate this, start out by assuming that the columns of $X$ are orthonormal, and apply the Householder algorithm to this matrix. That might seem to be a strange thing to do, because the usual point of the Householder algorithm is to obtain an orthonormal basis for the columns of $X$. But that is not the goal here: the interest is in the orthocomplement. Notice that since $H_p ... H_2 H_1$ is an isometry and therefore preserves orthogonality, $H_p ... H_2 H_1 X$ has orthonormal columns if $X$ does, so the upper triangular matrix of the algorithm in this case must be diagonal with $\pm 1$'s on the diagonal. Since $H_p ... H_2 H_1 X = \begin{bmatrix} T \\ 0 \end{bmatrix}$ is equivalent to

$H_p ... H_2 H_1 X T^{-1} = \begin{bmatrix} I_p \\ 0 \end{bmatrix}$, nothing is lost by making the assumptions of the following lemma.

**Lemma 10:** *Let $X$ be an $n \times p$ matrix with orthonormal columns, with $p \leq n$. Let $H_1, ..., H_p$ be the sequence of reflections from the Householder factorization algorithm in the notation described above, such that $H_k ... H_1 X = \begin{bmatrix} I_k \\ 0 \end{bmatrix}$ for $k \leq p$ (this amounts to choosing the "to" vector be $e_k$ at step $k$, rather than $-e_k$). Then there are $k \times k$ matrices $S_k$, which only depend on the first $k$ rows and columns of $X$, such that*

*(1)* $H_k ... H_1 = I - \left( [x_1, ..., x_k] - \begin{bmatrix} I_k \\ 0 \end{bmatrix} \right) S_k \left( [x_1, ..., x_k] - \begin{bmatrix} I_k \\ 0 \end{bmatrix} \right)^T, 1 \leq k \leq p$.

*(2)* $\left( I_k - [x_1, ..., x_k]^{(k)} \right) S_k \left( I_k - [x_1, ..., x_k]^{(k)} \right) = I_k - [x_1, ..., x_k]^{(k)}, 1 \leq k \leq p$, *so*

$rank\, S_k \geq rank \left( I_k - [x_1, ..., x_k]^{(k)} \right)$. *If* $I_k - [x_1, ..., x_k]^{(k)}$ *is invertible,* $S_k = \left( I_k - [x_1, ..., x_k]^{(k)} \right)^{-1}$.

*(3) For $1 \leq k \leq p$, if $x \in \left( col\, [x_1, ..., x_k] \right)^\perp$, then*

$\qquad$ *(i)* $\left( I_k - [x_1, ..., x_k]^{(k)} \right) S_k x^{(k)} = x^{(k)}$, *and*

$\qquad$ *(ii)* $\left( H_k ... H_1 \right)_{(k)} x = x_{(k)} + [x_1, ..., x_k]_{(k)} S_k x^{(k)}$.

*This second statement, when $k = p$, is what is used in the statistical application.*

*(4)* $S_k^T \left( I_k - [x_1, ..., x_k]^{(k)} \right)^T S_k = S_k$ *for* $1 \leq k \leq p$;

$\qquad$ *at every step, the rank of $S_k$ equals the rank of $I_k - [x_1, ..., x_k]^{(k)}$.*

*(5) The $S_k$'s satisfy the recursion*

$$S_{k+1} = \begin{bmatrix} S_k & 0 \\ 0 & 0 \end{bmatrix} + \frac{1}{1 - x_{k+1,k+1} - x_{k+1,1:k} S_k x_{k+1}^{(k)}} \begin{bmatrix} S_k x_{k+1}^{(k)} \\ 1 \end{bmatrix} \begin{bmatrix} x_{k+1,1:k} S_k & 1 \end{bmatrix}$$

*if $v_{k+1} \neq 0$ (equivalent to $1 - x_{k+1,k+1} - x_{k+1,1:k} S_k x_{k+1}^{(k)} \neq 0$ ), else $S_{k+1} = \begin{bmatrix} S_k & 0 \\ 0 & 0 \end{bmatrix}$. The rank of $S_{k+1}$ is one more than the rank of $S_k$ if and only if $v_{k+1} \neq 0$. $Rank\left( S_k \right) = \#\left\{ i \leq k : v_i \neq 0 \right\} = \#\left\{ i \leq k : H_i \neq I_n \right\}$.*

In (5), $x_{k+1,1:k}$ denotes the row vector consisting of the $(k+1)^{st}$ row of $X$ and the first $k$ columns.

*Remark.* There is a well-known formula for the inverse of a block matrix (see the Wikipedia entry for block matrix), which for diagonal blocks being $k \times k$ and $1 \times 1$ can be written as

$$\begin{bmatrix} A & b \\ c^T & d \end{bmatrix}^{-1} = \begin{bmatrix} A^{-1} & 0 \\ 0 & 0 \end{bmatrix} + \frac{1}{d - c^T A^{-1} b} \begin{bmatrix} -A^{-1} b \\ 1 \end{bmatrix} \begin{bmatrix} -c^T A^{-1} & 1 \end{bmatrix}$$

with the second term having rank one. Write $I_{k+1} - X^{(k+1)}$ in block matrix form as

$$I_{k+1} - X^{(k+1)} = \begin{bmatrix} I_k - [x_1, ... x_k]^{(k)} & -x_{k+1}^{(k)} \\ -x_{k+1,1:k} & 1 - x_{k+1,k+1} \end{bmatrix} = \begin{bmatrix} A & b \\ c^T & d \end{bmatrix} \text{ where}$$

$A = I_k - [x_1, ... x_k]^{(k)}, \ b = -x_{k+1}^{(k)} \in \mathbb{R}^k, c = -x_{k+1,1:k}^{\ T} \in \mathbb{R}^k, d = 1 - x_{k+1,k+1} \in \mathbb{R}$.

A glance at this formula and the recursion (5) shows that if $S_k = \left( I_k - [x_1, ... x_k]^{(k)} \right)^{-1} = A^{-1}$, and if

$v_{k+1} \neq 0$, then $S_{k+1} = \left( I_{k+1} - [x_1, ... x_{k+1}]^{(k+1)} \right)^{-1}$ (and that $I_{k+1} - [x_1, ... x_{k+1}]^{(k+1)}$ is invertible). Conversely, if $I_{k+1} - [x_1, ... x_{k+1}]^{(k+1)}$ is invertible, then the block inversion formula shows that $d - c^T A^{-1} b$ is not zero, so $v_{k+1} \neq 0$. In the non-singular case, our formula amounts to computing the inverse recursively by using that block matrix inversion formula, ending with $S_p = \left( I_p - [x_1, ... x_p]^{(p)} \right)^{-1}$. But the recursion still works to compute the $S_k$'s even in the singular case. It represents the result built up as a sum of rank-one operators.

Before proving Lemma 10, it will be used to prove Theorems 1 and 2.

*Proof of Theorem* 2. Suppose $XC^{-1}$ has orthonormal columns. Applying Lemma 10 to this matrix, there is a matrix $S_p{}'$ having the same rank as $I_p - X^{(p)} C^{-1}$, such that from part (3)(ii) of the lemma,

$$\left( H_p ... H_1 \right)_{(p)} x = x_{(p)} + X_{(p)} C^{-1} S_p{}' x^{(p)} \text{ for all } x \in \left( col\left( XC^{-1} \right) \right)^{\perp}.$$ But the column space of $XC^{-1}$ is the

same as the column space of $X$, and the rank of $C^{-1} S_p{}'$ is the same as the rank of $S_p{}'$, and the rank of

$I_p - X^{(p)} C^{-1}$ is the same as the rank of $C - X^{(p)}$. Letting $M^T = \left( H_p ... H_1 \right)_{(p)}$ and $S = C^{-1} S_p{}'$, the first statement of Theorem 2 is obtained. If $C - X^{(p)}$ is non-singular, then $I_p - X^{(p)} C^{-1}$ is non-singular,

and $S_p{}' = \left( I_p - X^{(p)} C^{-1} \right)^{-1}$ by (2) of Lemma 10. But then

$$S = C^{-1} S_p{}' = C^{-1} \left( I_p - X^{(p)} C^{-1} \right)^{-1} = \left( C - X^{(p)} \right)^{-1},$$ which proves the second statement of Theorem 2.

It is clear from (5) of Lemma 10 that $S$ can be computed in $O(p^3)$ time by that recursive algorithm; in the non-singular case, that is also the time needed to compute the inverse matrix. Just one more simple lemma is needed:

**Lemma 11.** If $M$ is any $p \times p$ matrix, there is a diagonal matrix $D = diag\left( d_1, d_2, ..., d_p \right)$ where

$d_i = \pm 1$ for each $i$, such that $D - M$ is non-singular. $D$ can be found in $O\left( p^3 \right)$ time.

*Proof.* Assume $d_1, d_2, ..., d_k$ have been found so that $A = diag\left( d_1, d_2, ..., d_k \right) - m_{1:k,1:k}$, the $k \times k$

principal submatrix of $D - M$, is non-singular, and that $A^{-1}$ is found. The block matrix inversion formula mentioned earlier shows that $diag\left( d_1, d_2, ..., d_{k+1} \right) - m_{1:k+1,1:k+1}$ is invertible, and how to compute its inverse in $O\left(k^2\right)$ time, if $d_{k+1} - m_{k+1,k+1} - c^T A^{-1} b \neq 0$, where $c^T = m_{k+1,1:k}$, $b = m_{1:k,k+1}$. For at least one choice of $d_{k+1} = \pm 1$, this must be true. By induction the result follows. □

Use this lemma to find $D$ such that $D - X^{(p)} C^{-1}$ is non-singular, and therefore $DC - X^{(p)}$ is non-singular, completing the proof of Theorem 2. □

*Proof of Theorem* 1. This actually follows from the proof of Theorem 2, with the observation that $H_p...H_1$ of Theorem 10 is the same as the $H_p...H_1$ of the Householder algorithm that produced $H_p...H_1 X = \begin{bmatrix} T \\ 0 \end{bmatrix}$, so the $M$ of Theorem 2 is the $U_2$ of Theorem 1; and $T - X^{(p)}$ is non-singular by Lemma 8. □

*Proof of Lemma 10.* This is straightforward, though long to write out because it has several parts. To prove (1), make the induction hypothesis

$$H_k...H_1 = I - \left(\left[x_1,...,x_k\right] - \begin{bmatrix} I_k \\ 0 \end{bmatrix}\right) S_k \left(\left[x_1,...,x_k\right] - \begin{bmatrix} I_k \\ 0 \end{bmatrix}\right)^T \text{ for some } k \times k \text{ matrix } S_k.$$

By definition in the Householder algorithm,

$$v_{k+1} = H_k...H_1 x_{k+1} - e_{k+1} = x_{k+1} - e_{k+1} - \left(\left[x_1,...,x_k\right] - \begin{bmatrix} I_k \\ 0 \end{bmatrix}\right) S_k \left(\left[x_1,...,x_k\right] - \begin{bmatrix} I_k \\ 0 \end{bmatrix}\right)^T x_{k+1}$$

$$= x_{k+1} - e_{k+1} - \left(\left[x_1,...,x_k\right] - \begin{bmatrix} I_k \\ 0 \end{bmatrix}\right) S_k \left(-x_{k+1}^{(k)}\right)$$

$$= \left(\left[x_1,...,x_{k+1}\right] - \begin{bmatrix} I_{k+1} \\ 0 \end{bmatrix}\right) q_{k+1}, \text{ where } q_{k+1} = \begin{bmatrix} S_k x_{k+1}^{(k)} \\ 1 \end{bmatrix}.$$

It is possible that $v_{k+1} = 0$, so $H_{k+1} = I$, and $S_{k+1} = \begin{bmatrix} S_k & 0 \\ 0 & 0 \end{bmatrix}$. If not,

$$H_{k+1}H_k...H_1 = \left( I - \frac{2}{\|v_{k+1}\|^2} v_{k+1} v_{k+1}{}^T \right)\left( I - \left( [x_1,...,x_k] - \begin{bmatrix} I_k \\ 0 \end{bmatrix} \right) S_k \left( [x_1,...,x_k] - \begin{bmatrix} I_k \\ 0 \end{bmatrix} \right)^T \right)$$

$$= I - \left( [x_1,...,x_{k+1}] - \begin{bmatrix} I_{k+1} \\ 0 \end{bmatrix} \right)\begin{bmatrix} S_k & 0 \\ 0 & 0 \end{bmatrix}\left( [x_1,...,x_{k+1}] - \begin{bmatrix} I_{k+1} \\ 0 \end{bmatrix} \right)^T$$

$$- \frac{2}{\|v_{k+1}\|^2}\left( [x_1,...,x_{k+1}] - \begin{bmatrix} I_{k+1} \\ 0 \end{bmatrix} \right) q_{k+1} q_{k+1}{}^T \left( [x_1,...,x_{k+1}] - \begin{bmatrix} I_{k+1} \\ 0 \end{bmatrix} \right)^T$$

$$\times \left\{ I - \left( [x_1,...,x_{k+1}] - \begin{bmatrix} I_{k+1} \\ 0 \end{bmatrix} \right)\begin{bmatrix} S_k & 0 \\ 0 & 0 \end{bmatrix}\left( [x_1,...,x_{k+1}] - \begin{bmatrix} I_{k+1} \\ 0 \end{bmatrix} \right)^T \right\}$$

$$= I - \left( [x_1,...,x_{k+1}] - \begin{bmatrix} I_{k+1} \\ 0 \end{bmatrix} \right) S_{k+1} \left( [x_1,...,x_{k+1}] - \begin{bmatrix} I_{k+1} \\ 0 \end{bmatrix} \right)^T,$$

where

$$S_{k+1} = \begin{bmatrix} S_k & 0 \\ 0 & 0 \end{bmatrix} + \frac{2}{\|v_{k+1}\|^2} q_{k+1} q_{k+1}{}^T \left\{ I_{k+1} - \left( [x_1,...,x_{k+1}] - \begin{bmatrix} I_{k+1} \\ 0 \end{bmatrix} \right)^T \left( [x_1,...,x_{k+1}] - \begin{bmatrix} I_{k+1} \\ 0 \end{bmatrix} \right)\begin{bmatrix} S_k & 0 \\ 0 & 0 \end{bmatrix} \right\}$$

$$= \begin{bmatrix} S_k & 0 \\ 0 & 0 \end{bmatrix} + \frac{2}{\|v_{k+1}\|^2} q_{k+1} q_{k+1}{}^T \left\{ I_{k+1} - \left( I_{k+1} - [x_1,...,x_{k+1}]^{(k+1)T} + I_{k+1} - [x_1,...,x_{k+1}]^{(k+1)} \right)\begin{bmatrix} S_k & 0 \\ 0 & 0 \end{bmatrix} \right\}.$$

That proves the induction step, and part (1) of the lemma. The last line shows that $S_{k+1}$ only depends on the first $k+1$ rows and columns of $X$. It also gives a recursion computing $S_{k+1}$ from $S_k$ in $O\left((k+1)^2\right)$ steps, so $S_p$ is computable in $O\left(p^3\right)$ steps. The recursion can be simplified to the statement in (5), which will be done later.

Next,

$$H_k...H_1[x_1,...,x_k] = [x_1,...,x_k] - \left( [x_1,...,x_k] - \begin{bmatrix} I_k \\ 0 \end{bmatrix} \right) S_k \left( [x_1,...,x_k] - \begin{bmatrix} I_k \\ 0 \end{bmatrix} \right)^T [x_1,...,x_k]$$

$$= [x_1,...,x_k] - \left( [x_1,...,x_k] - \begin{bmatrix} I_k \\ 0 \end{bmatrix} \right) S_k \left( I_k - [x_1,...,x_k]^{(k)} \right).$$

But also

$$H_k...H_1[x_1,...,x_k] = \begin{bmatrix} I_k \\ 0 \end{bmatrix}, \text{ so equating the first k rows of the right sides above,}$$

$$I_k = [x_1,...,x_k]^{(k)} - \left( [x_1,...,x_k]^{(k)} - I_k \right) S_k \left( I_k - [x_1,...,x_k]^{(k)} \right),$$

which when rearranged proves (2) of the lemma.

Next, the QR factorization $\left[x_1,...,x_k\right] = \left(H_k...H_1\right)^T \begin{bmatrix} I_k \\ 0 \end{bmatrix}$ implies $H_k...H_1 = \begin{bmatrix} \left[x_1,...,x_k\right]^T \\ U_k^T \end{bmatrix}$, where the columns of $U_k$ are an orthonormal basis for $\left(col\left[x_1,...,x_k\right]\right)^\perp$. Let $x \in \left(col\left[x_1,...,x_k\right]\right)^\perp$.

$$H_k...H_1 x = \begin{bmatrix} \left[x_1,...,x_k\right]^T \\ U_k^T \end{bmatrix} x = \begin{bmatrix} 0 \\ U_k^T x \end{bmatrix} = x - \left(\left[x_1,...,x_k\right] - \begin{bmatrix} I_k \\ 0 \end{bmatrix}\right) S_k \left(\left[x_1,...,x_k\right] - \begin{bmatrix} I_k \\ 0 \end{bmatrix}\right)^T x$$

$$= x + \left(\left[x_1,...,x_k\right] - \begin{bmatrix} I_k \\ 0 \end{bmatrix}\right) S_k x^{(k)}.$$

The first $k$ rows of this equation give

$0 = x^{(k)} + \left(\left[x_1,...,x_k\right]^{(k)} - I_k\right) S_k x^{(k)}$, or $\left(I_k - \left[x_1,...,x_k\right]^{(k)}\right) S_k x^{(k)} = x^{(k)}$, which proves (3)(i) of the

lemma. The last $n-k$ rows give $U_k^T x = x_{(k)} + \left[x_1,...,x_k\right]_{(k)} S_k x^{(k)}$, which is (3)(ii).

*Remark*: At this point, (1) through (3) of the lemma is sufficient for proving Theorem 1, and Theorem 2 in the non-singular case, because (2) gives the simple formula for $S_p$ when $I_p - X^{(p)}$ is invertible. The rest of the proof of Lemma 10 is straightforward but somewhat long and tedious to write out. We leave the proof to the reader if interested in the singular case. □

## 3. Application to Statistics.

Assume the regression setup given in the introduction. In classical statistics, the chi-square distribution of the sum of squared residuals is often proved using Cochran's theorem (see e.g. Hogg, McKean and Craig [**2**, p. 520]), which concerns the distribution of quadratic forms in independent standard normal variables. There is the standard partition of the sums of squares,

$Q = \left(Y - X\beta\right)^T\left(Y - X\beta\right) = \left(\hat\beta - \beta\right)^T\left(X^T X\right)\left(\hat\beta - \beta\right) + \left(Y - X\hat\beta\right)^T\left(Y - X\hat\beta\right) = Q_1 + Q_2$. The matrix

of $Q$ as a quadratic form in $Y - X\beta$ is $I_n$. One sees that

$\left(\hat\beta - \beta\right)^T\left(X^T X\right)\left(\hat\beta - \beta\right) = \left(Y - X\beta\right)^T X\left(X^T X\right)^{-1} X^T\left(Y - X\beta\right)$, so the matrix of $Q_1$ as a quadratic

form in $Y - X\beta$ is $B_1 = X\left(X^T X\right)^{-1} X^T$, which has rank $p$. $I_n = B_1 + B_2$, where $B_2$ is the matrix of the

$Q_2 = R^T R$ as a quadratic form in $Y - X\beta$. Now $B_1 X = X\left(X^T X\right)^{-1} X^T X = X$, and $I_n X = X$, so

$B_2 X = 0$, showing $rank(B_2) \leq n - p$, so in fact $rank(B_2) = n - p$. Since $rank(B_1) + rank(B_2) = n$, it

follows from Cochran's Theorem that $Q_2 / \sigma^2$ has the chi-square distribution with $n - p$ degrees of freedom (and $Q_1 / \sigma^2$ has the chi-square distribution with $p$ degrees of freedom, and $Q_1$ and $Q_2$ are

independent). This is the classical theory, which arrives at this conclusion by showing the existence of

some $n \times (n - p)$ matrix $M$ with orthonormal columns such that $B_2 = MM^T$, so that

$W = M^T\left(Y - X\beta\right) \sim N(0, \sigma^2 I_{n-p})$ and $W^T W = Q_2 = R^T R$. On the left is the sum of squares of

$n - p$ independent terms, and on the right the sum of squares of the $n$ residuals, which are however not independent.

It is worth noting that since $B_2 X = 0$ implies $M^T X = 0$, $W = M^T Y$, with no appearance of the unknown parameters $\beta$, and also $W = M^T \left( Y - X \hat{\beta} \right) = M^T R$, also not involving unknown parameters $\beta$. This means that $W$ is actually a statistic, computable from the data without knowing the parameters. But the classical texts are only focused on the distribution of $W^T W$, which does not require actually finding a way to compute $W$.

Yiping Cheng [**1**], motivated by pedagogy, exhibits such a representation for the case $p = 1$ (Student's Theorem). He gives a formula for an $n \times n$ orthogonal matrix whose last $n - 1$ rows are what we are calling $M^T$: $M$ is an $n \times (n-1)$ matrix with orthonormal columns which are orthogonal to the vector of ones, such that if $W = M^T \left( Y - \mathbf{1} \mu \right)$, then $\sum_{j=1}^{n-1} W_j^2 = \sum_{j=1}^{n} \left( Y_j - \overline{Y} \right)^2$. Because of the orthonormal columns, $W \sim N(0, \sigma^2 I_{n-1})$. Since $M^T \mathbf{1} = 0$, actually $W = M^T Y$ and is computable without knowing $\mu$, so it is a computable representation.

Cheng creates his matrix by induction and insight, focusing on keeping the matrix orthogonal. But actually, *any* factorization of $B_2$ as $B_2 = MM^T$, with $M$ having $n - 1$ columns, forces $M$ to have orthonormal columns! This is a consequence of $B_2$ being idempotent, as will be shown below in Lemma 12. The fact that $B_2$ is idempotent follows from this simple lemma, which is essentially in the material in Rao [**3**, p. 28], or can easily be proved by the reader:

*Lemma*. Let $B$ be a symmetric $n \times n$ matrix. Then $rank(B) + rank(I_n - B) = n$ if and only if $B$ is idempotent ( meaning $B^2 = B$; equivalently, the eigenvalues of $B$ are zero or one).

The following simple lemma gives an easy way of obtaining Cheng's matrix. Perhaps this rigidity in the factorization of an idempotent symmetric matrix is of other interest.

**Lemma 12**. *Let $B$ be an idempotent, symmetric $n \times n$ matrix of rank $r$. If $M$ is any matrix with $r$ columns such that $B = MM^T$, then $M$ has orthonormal columns.*

*Proof.* $B = MM^T$, so $MM^T = B = B^2 = MM^T MM^T$, and then $M^T MM^T M = M^T MM^T MM^T M$. But $M$ is of rank $r$, so $M^T M$ is an $r \times r$ matrix of rank $r$, and thus invertible. Cancelling leaves $M^T M = I_r$, which proves the lemma. □

To obtain Cheng's matrix, simple Gaussian elimination can be used to obtain the $B_2 = LDL^T$ factorization as taught in [**4**, p. 51], with diagonal matrix $D$, and set $M = LD^{1/2}$. Nothing has to be done to create orthonormal columns: they will just come out that way automatically, by Lemma 12. Here are the steps carried out for Student's problem.

For $p = 1$, $B_1$ is rank 1, and is simply $B_1 = \dfrac{1}{n}\mathbf{1}\mathbf{1}^{\mathrm{T}}$, where $\mathbf{1}^{\mathrm{T}} = \begin{bmatrix} 1 & 1 & \dots & 1 \end{bmatrix}$ is the vector of $n$ ones. So

$$B_2 = I_n - B_1 = \begin{bmatrix} \dfrac{n-1}{n} & -\dfrac{1}{n} & \dots & -\dfrac{1}{n} \\ -\dfrac{1}{n} & \dfrac{n-1}{n} & \dots & -\dfrac{1}{n} \\ \dots & \dots & \dots & \dots \\ -\dfrac{1}{n} & -\dfrac{1}{n} & \dots & \dfrac{n-1}{n} \end{bmatrix}.$$

One step of row reduction using upper left corner as pivot gives

$$B_2 \sim \begin{bmatrix} \dfrac{n-1}{n} & -\dfrac{1}{n} & \dots & -\dfrac{1}{n} \\ 0 & \dfrac{n-2}{n-1} & \dots & -\dfrac{1}{n-1} \\ \dots & \dots & \dots & \dots \\ 0 & -\dfrac{1}{n-1} & \dots & \dfrac{n-2}{n-1} \end{bmatrix}; \quad B_2 \sim U = \begin{bmatrix} \dfrac{n-1}{n} & -\dfrac{1}{n} & \dots & -\dfrac{1}{n} & -\dfrac{1}{n} \\ 0 & \dfrac{n-2}{n-1} & \dots & -\dfrac{1}{n-1} & -\dfrac{1}{n-1} \\ \dots & \dots & \dots & \dots & \dots \\ 0 & 0 & \dots & \dfrac{1}{2} & -\dfrac{1}{2} \\ 0 & 0 & \dots & 0 & 0 \end{bmatrix}$$ after n-1 steps.

Factoring out the diagonal terms of $U$ from each row and discarding the zero row makes this just the transpose of the lower triangular $L$ factor in the "rank revealing" factorization $B_2 = LDL^T$, where

$$D = \begin{bmatrix} \dfrac{n-1}{n} & 0 & \dots & 0 \\ 0 & \dfrac{n-2}{n-1} & \dots & 0 \\ \dots & \dots & \dots & \dots \\ 0 & 0 & \dots & \dfrac{1}{2} \end{bmatrix} \text{ is } (n-1)\times(n-1), \text{ and } L^T = \begin{bmatrix} 1 & -\dfrac{1}{n-1} & \dots & -\dfrac{1}{n-1} & -\dfrac{1}{n-1} \\ 0 & 1 & \dots & -\dfrac{1}{n-2} & -\dfrac{1}{n-2} \\ \dots & \dots & \dots & \dots & \dots \\ 0 & 0 & \dots & 1 & -1 \end{bmatrix},$$

where $L$ is $n \times (n-1)$. Let $M = LD^{1/2}$. Then $B_2 = MM^T$. Note that the columns of $L$ are orthogonal, and the columns of $M$ are orthonormal. Lemma 12 says this *had* to happen.

This is the same as Cheng's solution, which satisfies his goal of achieving a concrete representation. But the calculation of $M^T Y$ using this matrix takes $O(n^2)$ operations, so is not particularly simple.

**A better way with Householder.** Though Lemma 12 made it fairly easy to obtain $n-1$ orthonormal vectors orthogonal to the vector $\dfrac{1}{\sqrt{n}}\mathbf{1}$, it still took $n-1$ row operations, and this is the wrong way to do it: that is the insight mentioned in the introduction. The last $n-1$ columns of a single reflection matrix will suffice to give an orthonormal basis for the orthocomplement of $span\left\{\dfrac{1}{\sqrt{n}}\mathbf{1}\right\}$. Let $H = I - 2\dfrac{vv^T}{\|v\|^2}$, where $v = \dfrac{1}{\sqrt{n}}\mathbf{1} - e_1$. Then $H\left(\dfrac{1}{\sqrt{n}}\mathbf{1}\right) = e_1$, so $\dfrac{1}{\sqrt{n}}\mathbf{1} = H^T e_1 = H e_1$ which is the first column of $H$, so $H = \left[\begin{array}{cc}\dfrac{1}{\sqrt{n}}\mathbf{1} & M\end{array}\right]$ is an orthogonal matrix, and $M$ is an $n \times (n-1)$ matrix with orthonormal columns which are orthogonal to $\dfrac{1}{\sqrt{n}}\mathbf{1}$. Compute $\|v\|^2 = 2\left(1 - \dfrac{1}{\sqrt{n}}\right)$, so

$H = I - 2\dfrac{vv^T}{\|v\|^2} = I - \dfrac{\sqrt{n}}{\sqrt{n}-1}\left(\dfrac{1}{\sqrt{n}}\mathbf{1} - e_1\right)\left(\dfrac{1}{\sqrt{n}}\mathbf{1} - e_1\right)^T$, and just keep the last $n-1$ columns of this to

get $M$. This is a very simple matrix, the identity minus a simple rank-one matrix. But it is better than that: the only interest is in applying $M^T$ to the vector $R$ of residuals, to get $W$. Since $R$ is perpendicular to $\dfrac{1}{\sqrt{n}}\mathbf{1}$, $H^T R = R + \dfrac{\sqrt{n}}{\sqrt{n}-1}\left(\dfrac{1}{\sqrt{n}}\mathbf{1} - e_1\right)R_1$, and keeping the last $n-1$ rows,

$W = R_{(1)} + \dfrac{R_1}{\sqrt{n}-1}\mathbf{1}_{(1)}$. In coordinates, $W_j = R_{j+1} + \dfrac{R_1}{\sqrt{n}-1}$, $j = 1,\ldots,n-1$.

Motivated by how simple this turned out, the same idea was carried out for $p = 2$, "slope-intercept" straight-line regression, with $X$ made to have orthonormal columns before starting, working out the product of the two Householder reflections and applying it to $R$. Again, it collapsed into something simpler than expected: merely a correction term to $R$, depending only on the first two components of $R$, and the elements of the inverse matrix of $I - X^{(2)}$ were recognized in the formula. That motivated the general results of Lemma 10 and Theorems 1 and 2, to be applied to the multivariate regression problem in Theorem 3. In the Introduction, it was done the other way around; Theorem 3 was used to get a solution for Student's problem, and for the case $p = 2$. But that's not the order in which the discovery took place.

*Proof of Theorem 3.* Let $M$ be an $n \times (n-p)$ matrix whose columns are an orthonormal basis for

$col(X)^\perp$, and let $K$ be an $n \times p$ matrix whose columns are an orthonormal basis for $col(X)$, so

$U = \left[\begin{array}{cc}K & M\end{array}\right]$ is an orthogonal matrix. Let $W = M^T R$. The columns of $K$ are perpendicular to $R$, so

$R^T R = R^T U U^T R = \left[\begin{array}{cc}R^T K & R^T M\end{array}\right]\left[\begin{array}{c}K^T R \\ M^T R\end{array}\right] = \left[\begin{array}{cc}0 & W^T\end{array}\right]\left[\begin{array}{c}0 \\ W\end{array}\right] = W^T W$. Since the columns of $M$ are

perpendicular to the columns of $X$, $W = M^T R = M^T\left(Y - X\hat{\beta}\right) = M^T Y = M^T\left(Y - X\beta\right)$, so

$E\left[WW^T\right] = E\left[M^T\left(Y - X\beta\right)\left(Y - X\beta\right)^T M\right] = M^T \sigma^2 I_n M = \sigma^2 M^T M = \sigma^2 I_{n-p}$.

Theorem 3 now follows immediately from Theorems 1 and 2. □

*Proof of Corollary 5.* Since $X$ has orthonormal columns, if Householder is done not worrying about signs and just letting "to" be $e_1$ and $e_2$, then $T = I_2$. If $I_2 - X^{(2)}$ is non-singular, then

$S = \left(I_2 - X^{(2)}\right)^{-1}$; it is not necessary to actually carry out any Householder steps to get the answer.

$$S = \left(I_2 - X^{(2)}\right)^{-1} = \begin{bmatrix} 1 - \dfrac{1}{\sqrt{n}} & -t_1 \\[2mm] -\dfrac{1}{\sqrt{n}} & 1 - t_2 \end{bmatrix}^{-1} = \frac{1}{\det\left(I_2 - X^{(2)}\right)} \begin{bmatrix} 1 - t_2 & t_1 \\[2mm] \dfrac{1}{\sqrt{n}} & 1 - \dfrac{1}{\sqrt{n}} \end{bmatrix}$$

$$= \frac{1}{\left(\sqrt{n} - 1\right)\left(1 - t_2\right) - t_1} \begin{bmatrix} \sqrt{n}\left(1 - t_2\right) & \sqrt{n}\,t_1 \\[2mm] 1 & \sqrt{n} - 1 \end{bmatrix},$$

so

$XSR^{(2)} = \dfrac{1}{\left(\sqrt{n} - 1\right)\left(1 - t_2\right) - t_1}\left\{\left[\left(1 - t_2\right)R_1 + t_1 R_2\right]\mathbb{1} + \left[R_1 + \left(\sqrt{n} - 1\right)R_2\right]t\right\}$, which gives (a) for the

non-singular case. However, $I_2 - X^{(2)}$ could be singular, which occurs precisely when

$\left(\sqrt{n} - 1\right)\left(1 - t_2\right) - t_1 = 0$. This can happen in only one way:

$t_1 = \dfrac{1}{\sqrt{n}}, t_2 = 1 - \dfrac{1}{\sqrt{n}\left(\sqrt{n} - 1\right)}, t_i = -\dfrac{1}{\sqrt{n}\left(\sqrt{n} - 1\right)}, i = 3, \ldots, n$. This may be shown using the observation from Lemma 8 that $I_2 - X^{(2)}$ is singular if and only if $v_2$ is zero; the details are omitted.

The recursion formula from Lemma 10 gives

$S_1 = \left(1 - x_{11}\right)^{-1} = \left(1 - \dfrac{1}{\sqrt{n}}\right)^{-1} = \dfrac{\sqrt{n}}{\sqrt{n} - 1}, S = \begin{bmatrix} S_1 & 0 \\ 0 & 0 \end{bmatrix} = \dfrac{\sqrt{n}}{\sqrt{n} - 1}\begin{bmatrix} 1 & 0 \\ 0 & 0 \end{bmatrix}, XSR^{(2)} = \dfrac{R_1}{\sqrt{n} - 1}$, which gives

the singular case for (a). Unfortunately, when $I_2 - X^{(2)}$ is "almost" singular, one finds that the solution blows up and does not go gracefully to the singular case, so this is not so good.

Instead, one may can carry out the steps of the Householder algorithm with the standard sign choices, thereby avoiding a singular case. This gives

$v_1 = \dfrac{1}{\sqrt{n}}\mathbb{1} + e_1, \|v_1\|^2 = 2\left(1 + \dfrac{1}{\sqrt{n}}\right), H_1 x_2 = x_2 - \dfrac{2 v_1}{2\left(1 + \dfrac{1}{\sqrt{n}}\right)}t_1, \left(H_1 x_2\right)_2 = t_2 - \dfrac{1}{\sqrt{n} + 1}t_1,$

$s = \mathrm{sgn}(H_1 x_2)_2 = \mathrm{sgn}\left(\left(\sqrt{n} + 1\right)t_2 - t_1\right), v_2 = H_1 x_2 + s\,e_2$, so the "to" vectors are $-e_1$ and $-s\,e_2$, and

$T = \begin{bmatrix} -1 & 0 \\ 0 & -s \end{bmatrix}$. Then

$$S = \left(T - X^{(2)}\right)^{-1} = -\begin{bmatrix} 1 + \dfrac{1}{\sqrt{n}} & t_1 \\[2mm] \dfrac{1}{\sqrt{n}} & s + t_2 \end{bmatrix}^{-1} = \dfrac{-s}{\left(\sqrt{n}+1\right) + \left|\left(\sqrt{n}+1\right)t_2 - t_1\right|} \begin{bmatrix} \sqrt{n}\left(s + t_2\right) & -\sqrt{n}\,t_1 \\[2mm] -1 & \sqrt{n}+1 \end{bmatrix},$$

$$XSR^{(2)} = \dfrac{-s}{\left(\sqrt{n}+1\right) + \left|\left(\sqrt{n}+1\right)t_2 - t_1\right|} \left\{ \begin{array}{l} \left[\left(s + t_2\right)R_1 - t_1 R_2\right]\mathbf{1} \\[2mm] + \left[-R_1 + \left(\sqrt{n}+1\right)R_2\right]t \end{array} \right\}, \text{ which gives (b). } \square$$

## 4. Solution without Householder, by Hindsight

*Proof of Theorem 6*: It is elementary that $\mathrm{cov}(R_j, R_k) = -\dfrac{\sigma^2}{n}$ if $j \neq k$, and $Var(R_j) = \dfrac{n-1}{n}\sigma^2$. For

$j \neq k$, $\mathrm{cov}(W_j, W_k) = \mathrm{cov}(R_{j+1} + cR_1, R_{k+1} + cR_1) = -\dfrac{\sigma^2}{n} - 2c\dfrac{\sigma^2}{n} + c^2 \dfrac{n-1}{n}\sigma^2 = 0$ iff

$(n-1)c^2 - 2c - 1 = 0$. That's the condition for independence; they are clearly identically distributed

and normal, and for $c$ satisfying this equation, one finds $Var(W_j) = \sigma^2$. Compute

$$\sum_{j=1}^{n-1} W_j^2 = \sum_{j=1}^{n-1} \left(R_{j+1} + cR_1\right)^2 = \sum_{j=1}^{n-1} R_{j+1}^2 + 2cR_1 \sum_{j=1}^{n-1} R_{j+1} + (n-1)c^2 R_1^2$$

$$= \sum_{j=1}^{n-1} R_{j+1}^2 + 2cR_1(-R_1) + (n-1)c^2 R_1^2 = \sum_{j=1}^{n-1} R_{j+1}^2 + \left((n-1)c^2 - 2c\right)R_1^2,$$

which equals $\sum_{j=1}^{n} R_j^2$ for all $R$ iff $(n-1)c^2 - 2c = 1$. This is the same as the condition for independence, and the two solutions are $c = \dfrac{1}{\sqrt{n-1}}$ or $c = \dfrac{-1}{\sqrt{n+1}}$.

*Proof of Theorem 7*: Start with the condition $W^T W = R^T R$ to see what condition that imposes on $S$.
To simply the work, first assume that $X$ has orthonormal columns, so that

$X_{(p)}^T X_{(p)} + X^{(p)T} X^{(p)} = X^T X = I_p$. Also note $X^T R = 0 \Rightarrow X_{(p)}^T R_{(p)} = -X^{(p)T} R^{(p)}$. Then

$W^T W = \left(R_{(p)} + X_{(p)} SR^{(p)}\right)^T \left(R_{(p)} + X_{(p)} SR^{(p)}\right)$

$= R_{(p)}^T R_{(p)} + R_{(p)}^T X_{(p)} SR^{(p)} + R^{(p)T} S^T X_{(p)}^T R_{(p)} + R^{(p)T} S^T X_{(p)}^T X_{(p)} SR^{(p)}$

$= R_{(p)}^T R_{(p)} - R^{(p)T} X^{(p)} SR^{(p)} - R^{(p)T} S^T X^{(p)T} R^{(p)} + R^{(p)T} S^T \left(I_p - X^{(p)T} X^{(p)}\right) SR^{(p)}$

$= R_{(p)}^T R_{(p)} + R^{(p)T} \left(S^T \left(I_p - X^{(p)T} X^{(p)}\right)S - X^{(p)} S - S^T X^{(p)T}\right) R^{(p)}.$

A sufficient condition for this to equal $R^T R$ is

$$S^T \left(I_p - X^{(p)T} X^{(p)}\right)S - X^{(p)} S - S^T X^{(p)T} = I_p, \qquad (*)$$

and this is also necessary if it is to equal $R^T R$ for all possible $R$. Assuming $S$ is invertible, multiplying (*) on left and right by $S^{-T}$ and $S^{-1}$ resp.,

$I_p - X^{(p)T} X^{(p)} - S^{-T} X^{(p)} - X^{(p)T} S^{-1} = S^{-T} S^{-1}$, or $\left( S^{-1} + X^{(p)} \right)^T \left( S^{-1} + X^{(p)} \right) = I_p$. This is true iff $S^{-1} + X^{(p)} = Q$, where $Q$ is a $p \times p$ orthogonal matrix such that $Q - X^{(p)}$ is non-singular. Then $S = \left( Q - X^{(p)} \right)^{-1}$. The simplest is $S = \left( I_p - X^{(p)} \right)^{-1}$, if that is invertible.

Now consider the covariance condition. WLOG assume $\sigma^2 = 1$ to simplify the writing. The condition to be shown is $\operatorname{cov}(W, W) = E\left[ WW^T \right] = I_{n-p}$. The covariance of the residuals (see e.g. [3, p. 185)]) is $E\left[ RR^T \right] = I - X \left( X^T X \right)^{-1} X^T = I - XX^T$. In block matrices,

$E\left[ RR^T \right] = E\begin{bmatrix} R^{(p)} R^{(p)T} & R^{(p)} R_{(p)}{}^T \\ R_{(p)} R^{(p)T} & R_{(p)} R_{(p)}{}^T \end{bmatrix} = \begin{bmatrix} I_p - X^{(p)} X^{(p)T} & -X^{(p)} X_{(p)}{}^T \\ -X_{(p)} X^{(p)T} & I_{n-p} - X_{(p)} X_{(p)}{}^T \end{bmatrix}$. Thus

$E\left[ WW^T \right] = E\left[ \left( R_{(p)} + X_{(p)} SR^{(p)} \right) \left( R_{(p)} + X_{(p)} SR^{(p)} \right)^T \right]$

$= E\left[ R_{(p)} R_{(p)}{}^T + X_{(p)} SR^{(p)} R_{(p)}{}^T + R_{(p)} R^{(p)T} S^T X_{(p)}{}^T + X_{(p)} SR^{(p)} R^{(p)T} S^T X_{(p)}{}^T \right]$

$= I_{n-p} - X_{(p)} X_{(p)}{}^T - X_{(p)} SX^{(p)} X_{(p)}{}^T - X_{(p)} X^{(p)T} S^T X_{(p)}{}^T + X_{(p)} S\left( I_p - X^{(p)} X^{(p)T} \right) S^T X_{(p)}{}^T$

$= I_{n-p} - X_{(p)} \left( I_p + SX^{(p)} + X^{(p)T} S^T - SS^T + SX^{(p)} X^{(p)T} S^T \right) X_{(p)}{}^T$.

The requirement reduces to

$I_p + SX^{(p)} + X^{(p)T} S^T - SS^T + SX^{(p)} X^{(p)T} S^T = 0$.

Assuming invertibility, multiplying on the left and right by $S^{-1}$ and $S^{-T}$ resp., yields

$S^{-1} S^{-T} + X^{(p)} S^{-T} + SX^{(p)T} + X^{(p)} X^{(p)T} = I_p$, or $\left( S^{-1} + X^{(p)} \right) \left( S^{-1} + X^{(p)} \right)^T = I_p$.

This is actually equivalent to the condition we got previously when considering the sum of squares condition, even though the transposes are in reversed order, because a right inverse is also a left inverse. Both required conditions hold with $S = \left( Q - X^{(p)} \right)^{-1}$, if $Q$ is any $p \times p$ orthogonal matrix such that $Q - X^{(p)}$ is non-singular.

Finally, remove the requirement that $X$ has orthonormal columns, and suppose that $XC^{-1}$ has orthonormal columns. Replacing $X$ by $XC^{-1}$ in the formula just derived,

$W = R_{(p)} + X_{(p)} C^{-1} \left( Q - X^{(p)} C^{-1} \right)^{-1} R^{(p)} = R_{(p)} + X_{(p)} \left( QC - X^{(p)} \right)^{-1} R^{(p)}$

But $X(QC)^{-1} = \left( XC^{-1} \right) Q^T$ has orthonormal columns because $XC^{-1}$ does and $Q^T$ is orthogonal. Therefore, no increase in generality is obtained by putting $QC$ in place of just $C$ in the statement of Theorem 7. This completes the proof of Theorem 7. □